\documentclass{amsart}

\usepackage{amsmath,amssymb,amsthm,units, stmaryrd,stackrel,relsize,bm}

\usepackage{bussproofs}
\usepackage{proof}
\usepackage{yfonts,upgreek,units,stmaryrd,color,graphicx}
\usepackage[all]{xy}

\usepackage[
    colorlinks=true, 
    citecolor=red,
    linkcolor=blue,
    urlcolor=blue]{hyperref}
    
\newtheorem{theorem}{Theorem}[section]

\newtheorem{claim}[theorem]{Claim}

\newtheorem{proposition}[theorem]{Proposition}

\theoremstyle{definition}

\theoremstyle{remark}
\newtheorem{remark}[theorem]{Remark}

\DeclareSymbolFont{AMSb}{U}{msb}{m}{n}
\DeclareMathSymbol{\N}{\mathbin}{AMSb}{"4E}
\DeclareMathSymbol{\Z}{\mathbin}{AMSb}{"5A}
\DeclareMathSymbol{\R}{\mathbin}{AMSb}{"52}
\DeclareMathSymbol{\Q}{\mathbin}{AMSb}{"51}
\DeclareMathSymbol{\I}{\mathbin}{AMSb}{"49}
\DeclareMathSymbol{\C}{\mathbin}{AMSb}{"43}

\makeatletter
\def\Ddots{\mathinner{\mkern1mu\raise\p@
\vbox{\kern7\p@\hbox{.}}\mkern2mu
\raise4\p@\hbox{.}\mkern2mu\raise7\p@\hbox{.}\mkern1mu}}
\makeatother

\newcommand{\calc}{{\mathrm{LK}^{++}}}
\newcommand{\calcb}{{\mathrm{LK}^+}}
\newcommand{\calcint}{{\mathrm{LJ}^{+}}}
\newcommand{\ljcalc}{{\mathrm{LJ}^{++}}}
\newcommand{\lk}{{\mathrm{LK}}}
\newcommand{\lj}{{\mathrm{LJ}}}

\newcommand{\lkqs}{{\mathrm{LK}_{\texttt{shift}}}}
\newcommand{\lkep}{{\mathrm{LK}^\varepsilon}}

\newtheorem{exm}[theorem]{Example}
\newtheorem{cor}[theorem]{Corollary}

\theoremstyle{definition}
\newtheorem{defi}[theorem]{Definition}

\begin{document}
\title{Unsound Inferences Make Proofs Shorter}

\author{Juan P. Aguilera}
\address{Institute of Discrete Mathematics and Geometry,
Vienna University of Technology, 
Wiedner Hauptstra{\ss}e 8--10, 1040 Vienna, Austria
}

\author{Matthias Baaz}
\address{Institute of Discrete Mathematics and Geometry, 
Vienna University of Technology,
Wiedner Hauptstra{\ss}e 8--10, 1040 Vienna, Austria
}



\begin{abstract}
We give examples of calculi that extend Gentzen's sequent calculus $\lk$ by unsound quantifier inferences in such a way that (i) derivations lead only to true sequents, and (ii) proofs therein are non-elementarily shorter than $\lk$-proofs.
\end{abstract}

{\maketitle}
\section{Introduction} 
Consider the following argument:
\begin{enumerate}
\item \label{Int1} That Kurt G\"odel is Austrian entails that Kurt G\"odel is Austrian.
\item \label{Int2} Hence, that Kurt G\"odel is Austrian entails that everyone is Austrian.
\item \label{Int3} That is, if Kurt G\"odel is Austrian, then all people are Austrian.
\item \label{Int4} Therefore, there exists a person such that, if that person is Austrian, then all people are Austrian.
\end{enumerate}
The argument can be formalized as an instance of the following proof schema:
\begin{equation}  \label{ProofIntro}
\infer[]{\vdash \exists x\,\big( A(x) \to \forall y\, A(y)\big)}
{
	\infer[]{\vdash A(a) \to \forall y\, A(y)}
	{
		\infer[]{A(a) \vdash \forall y\, A(y)}{A(a) \vdash A(a)}
	}
}
\end{equation}
The study of formal proofs is motivated by our desire for deductive reasoning to be \textit{correct}, i.e., we wish that it be such that the procedures involved derive only true conclusions. The traditional way of ensuring this involves restricting proofs to those satisfying the two following properties:
\begin{itemize}
\item Inferences are \textit{sound},\footnote{Soundness is usually applied to derivations or logical systems. In this paper, we distinguish `soundness' from `correctness.' See below.} i.e., only true conclusions result from true premises.
\item Derivations are \textit{hereditary}, i.e., initial segments of proofs are proofs themselves.
\end{itemize}

In particular, one might impose certain \textit{characteristic-variable conditions} that restrict the circumstances under which quantifiers can be inferred.
We will use unsound rules based on weak characteristic-variable conditions to define enhanced (correct) logical calculi wherein derivations such as \eqref{ProofIntro} will be allowed. 

In our current treatment, the correctness of proofs can only be gauged by considering them in their entirety, as forfeiting soundness while maintaining correctness necessarily violates hereditariness. This is not unlike other aspects of reasoning that are traditionally not reflected in formal derivations. For example, hereditariness is not compatible with proofs by contradiction---one would not take as correct an initial segment of a proof obtained by interrupting it after an assumption directed towards a contradiction had been made, but before the contradiction had been reached.

To illustrate this, we note that the following subproof of the opening example should not be allowed as a proof in the system:
\begin{enumerate}
\item[1.] That Kurt G\"odel is Austrian entails that Kurt G\"odel is Austrian.
\item[2.] Hence, that Kurt G\"odel is Austrian entails that everyone is Austrian.
\end{enumerate}

The \textit{strong} inference from \ref{Int1} to \ref{Int2} is the only unsound one:

\begin{enumerate}
\item[1.] That Kurt G\"odel is Austrian entails that Kurt G\"odel is Austrian.
\item[3$'$] \label{Int3p} That is, if Kurt G\"odel is Austrian, then Kurt G\"odel is Austrian.
\item[4$'$] \label{Int4p} Therefore, there exists a person such that, if that person is Austrian, then Kurt G\"odel is Austrian.
\end{enumerate}

Here, the name `Kurt G\"odel' is simply manipulated syntactically, so its meaning plays no role. Moreover, the name does not even appear in sentence \ref{Int4}. It could well have been replaced by any other throughout the derivation with no effect---we call (the result of formalizing) this condition \textit{substitutability}. In our example, it is what allows us to mend the proof after a `lie' has been introduced in \ref{Int2}.

Formal proofs can often be taken without loss of generality to be \textit{regular}, i.e., any variable which is the \textit{characteristic variable}\footnote{See Section \ref{SectionQI} for preliminaries and definitions.} of a quantifier inference can be assumed to appear only before then. Forfeiting this condition allows us to accept the following proof:
\begin{align} \label{example2}
\infer[]{\vdash \exists x\, \big( A(x) \rightarrow A(f(x))\big)}{
\infer[]{\vdash A(a) \to A(f(a))}{
\infer[]{A(a) \vdash A(f(a))}{
	\infer[]{A(a) \vdash \forall x\, A(x)}{A(a) \vdash A(a)}
	&
	\infer[]{\forall x\, A(x) \vdash A(f(a))}{A(f(a)) \vdash A(f(a))}
}
}
}
\end{align}
In this example, $a$ is the characteristic variable of only one quantifier inference. We call this \textit{weak regularity}. Our main result is that a sequent calculus augmented with unsound quantifier inferences satisfying substitutability, weak regularity (or even a further weakening thereof), and a technical \textit{side-variable condition} is correct. This means that, although inferences therein may not be sound, the three (global) conditions are enough to guarantee that any lie introduced during the proof is eventually cleared.
Moreover, the calculus yields non-elementarily--shorter\footnote{Size can be measured, e.g., by counting the number of inferences in the proof.} cut-free proofs than standard proof systems. This is proved in Section \ref{SectSpeedup}. Some consequences of the speed-up theorem are explained in Section \ref{SectApplications}.

The three characteristic-variable conditions are motivated by the rules governing possible inferences in Hilbert's \textit{$\varepsilon$-calculus}, whose language contains no quantifiers (see Section \ref{SectEpsilon}). As a consequence, proofs therein are less restrictive. The weakened conditions thus result from attempting to pull back some extra freedom from the $\varepsilon$-calculus to first-order logic. They can also be thought of as a first attempt in tracking down and covering the spots where first-order logic proofs might break down without the usual quantifier-inference restrictions. In particular, it seems plausible that the conditions here can be further weakened, and the corresponding calculi can result even faster (see Section \ref{SectClosing}).

\section{Quantifier Inferences}\label{SectionQI}
We consider derivations in sequent calculi. Sequents are expressions of the form $\Gamma \vdash \Delta$, where $\Gamma$ and $\Delta$ represent collections of formulae. The sequent $\Gamma \vdash \Delta$ is usually interpreted as `if all of $\Gamma$ hold, then at least one of the formulae in $\Delta$ holds.' Equivalently, if $\Gamma = A_0, A_1, \hdots, A_n$ and $\Delta = B_0, B_1, \hdots, B_m$, then $\Gamma \vdash \Delta$ is interpreted as: 
\begin{equation}\label{seqdef}
\bigvee_{i\leq n} \lnot A_i \vee \bigvee_{i\leq m} B_i
\end{equation}

Inferences are expressions of the form 
\begin{equation*}
\infer[]{\Gamma' \vdash \Delta'}{\Gamma \vdash \Delta}
\end{equation*}
which are to be interpreted as `the sequent $\Gamma' \vdash \Delta'$ follows from the sequent $\Gamma \vdash \Delta$.' Derivations are trees whose leaves are axioms and whose non-leaf nodes are sequents obtained from their predecessors by inferences; they are to be interpreted as \textit{proofs} of the root sequent. In fact, we frequently refer to derivations as `proofs.'

If a formula is changed by an inference, we say that the formula is an \textit{auxiliary} formula of the inference. If so, then the resulting formula in the conclusion is called the \textit{critical} formula. An inference is a \textit{quantifier inference} if the (only) critical formula has a quantifier as its outermost logical symbol. Normally, quantifier inferences consist of substituting a formula\footnote{We will frequently use $Q$ to denote any unspecified quantifier.} $Qx\, A(x)$ for an instance thereof in a sequent, e.g., as in 
\begin{equation*}
\infer[]{\Gamma \vdash \Delta, \forall x\, A(x)}{\Gamma \vdash \Delta,  A(a)}
\end{equation*}
The \textit{polarity} of a quantifier in a sequent $\Gamma \vdash \Delta$ is defined as follows: rewrite $\Gamma \vdash \Delta$ in the form \eqref{seqdef} and then rewrite the resulting expression without implication symbols. A quantifier is \textit{positive} if it is universal (resp. existential) and under the scope of an even (resp. odd) number of negation symbols, and \textit{negative} otherwise.
A quantifier is \textit{strong} if it is positive and on the right-hand side of a sequent or negative and on the left-hand side of a sequent; it is \textit{weak} otherwise. 

If an inference yields a strongly-quantified formula $Qx\, A(x)$ from $A(a)$, where $a$ is a free variable, we say that $a$ is the \textit{characteristic variable} of the inference. We will denote free variables by letters $a, b, c, ...$ (in contrast, we denote bound variables by $x, y, z, ...$). We denote closed terms by $t, s, r, ...$ and variants thereof.

Let $\pi$ be any derivation. We say $b$ is a \textit{side variable} of $a$ in $\pi$ (written $a <_\pi b$) if $\pi$ contains a strong-quantifier inference of the form:
\begin{equation*}
\infer[]{\Gamma \vdash \Delta, \forall x\, A(x, b, \vec c)}{\Gamma \vdash \Delta,  A(a, b, \vec c)}
\end{equation*}
or of the form:
\begin{equation*}
\infer[]{\exists x\, A(x, b, \vec c), \Gamma \vdash \Delta}{A(a, b, \vec c), \Gamma \vdash \Delta}
\end{equation*}

The Skolemization of a first-order formula is defined by replacing every strongly quantified variable $y$ with a new function symbol $f_y(x_1, \hdots, x_n)$, where $x_1, \hdots, x_n$ are the weakly quantified variables such that $Q y$ appears in the scope of their quantifiers, and removing the quantifier $Qy$. We write $sk(A)$ for the Skolemization of $A$. Recall that whether a quantifier is weak or strong depends on the side of the sequent on which it is. Thus, $sk(A)$ means different things depending on whether $A$ appears in the antecedent or the succedent. The Skolemization of a sequent
$$A_0, \hdots, A_m \vdash A_{m+1}, \hdots, A_{m+n}$$
is defined as:
$$sk(A_0), \hdots, sk(A_m) \vdash sk(A_{m+1}), \hdots, sk(A_{m+n}).$$
In classical logic, a sequent is derivable if and only if its Skolemization is.

\begin{defi}[Suitable quantifier inference]
We say a quantifier inference is \emph{suitable for a proof $\pi$} if either it is a weak-quantifier inference, or the following three conditions are satisfied:
\begin{itemize}
\item (substitutability) the characteristic variable does not appear in the conclusion of $\pi$.
\item (side-variable condition) the relation $<_\pi$ is acyclic.
\item (weak regularity) the characteristic variable is not the characteristic variable of another strong-quantifier inference in $\pi$.
\end{itemize}
\end{defi}

We will work with various sequent calculi extending \textit{Propositional $\mathrm{LK}$}, the sequent calculus whose only axioms are $\bot \vdash$ and $A \vdash A$ (with $A$ atomic) and whose rules are:
\begin{align*}
\text{Structural rules:}\\
	\infer[WL]{\Gamma, A \vdash \Delta}{\Gamma \vdash \Delta} && 
	\infer[WR]{\Gamma \vdash \Delta, A}{\Gamma \vdash \Delta} \\
	\infer[CL]{\Gamma, A \vdash \Delta}{\Gamma, A, A \vdash \Delta} && 
	\infer[CR]{\Gamma \vdash \Delta, A}{\Gamma \vdash \Delta, A, A} \\	
	\infer[EL]{\Gamma_1, B, A, \Gamma_2 \vdash \Delta}{\Gamma_1, A, B, \Gamma_2 \vdash \Delta} && 
	\infer[ER]{\Gamma \vdash \Delta_1, B, A, \Delta_2}{\Gamma \vdash \Delta_1, A, B, \Delta_2} \\
	\infer[Cut]{\Gamma \vdash \Delta}{\Gamma \vdash \Delta, A & A, \Gamma \vdash \Delta} &&\\
\end{align*}
\begin{align*}
\text{Logical rules:}\\
	\infer[\wedge_1L]{\Gamma, A\wedge B \vdash \Delta}{\Gamma, A \vdash \Delta} && 
	\infer[\vee_1R]{\Gamma \vdash \Delta, A\vee B}{\Gamma \vdash \Delta, A} \\
	\infer[\wedge_2L]{\Gamma, A\wedge B \vdash \Delta}{\Gamma, B \vdash \Delta} && 
	\infer[\vee_2R]{\Gamma \vdash \Delta, A\vee B}{\Gamma \vdash \Delta, B} \\
	\infer[\vee L]{\Gamma, A\vee B \vdash \Delta}{\Gamma, A \vdash \Delta & \Gamma, B \vdash \Delta} && 
	\infer[\wedge R]{\Gamma \vdash \Delta, A\wedge B}{\Gamma \vdash \Delta, A & \Gamma \vdash \Delta, B} \\
	\infer[\to L]{\Gamma, A\to B \vdash \Delta}{\Gamma \vdash \Delta, A & \Gamma, B \vdash \Delta} && 
	\infer[\to R]{\Gamma \vdash \Delta, A\to B}{\Gamma, A \vdash \Delta, B}
\end{align*}
Negation $\lnot A$ is defined as $A \to \bot$. The names of the structural rules stand, respectively, for `weakening,' `contraction,' and `exchange.' First-order $\lk$ is the extension of Propositional LK obtained by adding quantifier inferences: 
\begin{align*}
	\infer[\forall L]{\Gamma, \forall x\, A(x) \vdash \Delta}{\Gamma, A(t) \vdash \Delta} && 
	\infer[\exists R]{\Gamma \vdash \Delta, \exists x\, A(x)}{\Gamma \vdash \Delta, A(t)} \\
	\infer[\exists L]{\Gamma, \exists x\, A(x) \vdash \Delta}{\Gamma, A(a) \vdash \Delta} && 
	\infer[\forall R]{\Gamma \vdash \Delta, \forall x\, A(x)}{\Gamma \vdash \Delta, A(a)}
\end{align*}
with the only restriction that in the inferences $\forall R$ and $\exists L$, the variable $a$ is not allowed to appear in the conclusion. In particular, note that every quantifier inference is suitable for every regular $\lk$-proof. Gentzen's famous \textit{Cut-Elimination Theorem} states that the cut rule is redundant in both propositional and first-order $\lk$.

\begin{defi}[$\calcb$]
The calculus $\calcb$ is defined like $\lk$, except that we instead allow all weak and strong quantifier inferences 
with the proviso that they be suitable for the proof.
\end{defi}

A further weakening of the characteristic-variable conditions gives rise to the notion of weak suitability:
\begin{defi}[Weakly suitable quantifier inference]
A quantifier inference is \textit{weakly suitable for a proof $\pi$} if either it is a weak-quantifier inference or it satisfies substitutability, the side-variable condition, and:
\begin{itemize}
\item (very weak regularity) whenever the characteristic variable is also the characteristic variable of another strong-quantifier inference in $\pi$, then it has the same critical formula.
\end{itemize}
\end{defi}

\begin{defi}[$\calc$]
The calculus $\calc$ is the extension of $\calcb$ that results from allowing all weakly suitable quantifier inferences.
\end{defi}

It is easy to find examples of sequents that are more easily provable in $\calcb$ than in $\lk$:

\begin{exm} \label{example25}
The sequent 
$$\forall x\, A(x) \rightarrow B \vdash \exists x\, (A(x) \rightarrow B)$$ is provable in $\lk$:
\begin{prooftree}
\def\fCenter{\mbox{\ $\vdash$\ }}
\Axiom$A(a) \fCenter A(a)$
\UnaryInf$A(a) \fCenter A(a), B$
\UnaryInf$\fCenter A(a), A(a) \rightarrow B$
\UnaryInf$\fCenter A(a), \exists x\, (A(x) \rightarrow B)$
\UnaryInf$\fCenter \exists x\, (A(x) \rightarrow B), A(a)$
\UnaryInf$\fCenter \exists x\, (A(x) \rightarrow B), \forall x\, A(x)$
\Axiom$B \fCenter B$
\BinaryInf$\forall x\, A(x) \rightarrow B \fCenter \exists x\, (A(x) \rightarrow B), B$
\UnaryInf$\forall x\, A(x) \rightarrow B, A(b) \fCenter \exists x\, (A(x) \rightarrow B), B$
\UnaryInf$\forall x\, A(x) \rightarrow B \fCenter \exists x\, (A(x) \rightarrow B), A(b) \rightarrow B$
\UnaryInf$\forall x\, A(x) \rightarrow B \fCenter \exists x\, (A(x) \rightarrow B), \exists x\, (A(x) \rightarrow B)$
\UnaryInf$\forall x\, A(x) \rightarrow B \fCenter \exists x\, (A(x) \rightarrow B)$
\end{prooftree}
However, one can find a shorter $\calcb$-proof:
\begin{prooftree}
\def\fCenter{\mbox{\ $\vdash$\ }}
\Axiom$A(a) \fCenter A(a)$
\UnaryInf$A(a) \fCenter \forall x\, A(x)$
\Axiom$B \fCenter B$
\BinaryInf$A(a), \forall x\, A(x) \rightarrow B \fCenter B$
\UnaryInf$\forall x\, A(x) \rightarrow B, A(a) \fCenter B$
\UnaryInf$\forall x\, A(x) \rightarrow B \fCenter A(a) \rightarrow B$
\UnaryInf$\forall x\, A(x) \rightarrow B \fCenter \exists x\, (A(x) \rightarrow B)$
\end{prooftree}
\qed
\end{exm}

Recall that a function on the natural numbers is \textit{elementary} if it can be defined by a quantifier-free formula from $+$, $\cdot$, and the function $x\mapsto 2^x$. 
By independent results of R. Statman \cite{statman79} (which were formalized in \cite{baazleitsch94}) and of V. P. Orevkov \cite{orevkov82}, the sizes of the smallest cut-free $\lk$-proofs of sequents of length $n$ are not bounded by any elementary function on $n$. Our main theorem is that cut-free $\calcb$-proofs are non-elementarily shorter than cut-free $\lk$-proofs:

\begin{theorem} \label{main}
There is no elementary function bounding the length of the shortest cut-free $\lk$-proof of a formula in terms of its shortest cut-free $\calcb$-proof.
\end{theorem}

An immediate consequence is the following:

\begin{cor} 
There is no elementary function bounding the length of the shortest cut-free $\lk$-proof of a formula in terms of its shortest cut-free $\calc$-proof.
\end{cor}
 
\noindent We prove Theorem \ref{main} in Section \ref{SectSpeedup}. First, we consider the question of the correctness of $\calc$. The proof of Theorem \ref{TheoremCorrectness} that we present here is due to the referee:

\begin{theorem} \label{TheoremCorrectness}
If a sequent is $\calc$-derivable, then it is already $\lk$-derivable.
\end{theorem}

\proof[Proof]
Let $\pi$ be an $\calc$-proof. Replace every unsound universal quantifier inference by a $\to L$ inference:
\begin{prooftree}
\def\fCenter{\mbox{\ $\vdash$\ }}
\Axiom$\Gamma \fCenter \Delta, A(a)$
\Axiom$\forall x\, A(x) \fCenter \forall x\, A(x)$
\BinaryInf$\Gamma, A(a) \to \forall x\, A(x) \fCenter \Delta, \forall x\, A(x)$
\end{prooftree}
Similarly replace every unsound existential quantifier by an $\to L$ inference
\begin{prooftree}
\def\fCenter{\mbox{\ $\vdash$\ }}
\Axiom$\exists x\, A(x) \fCenter \exists x\, A(x)$
\Axiom$A(a), \Gamma \fCenter \Delta$
\BinaryInf$\Gamma, \exists x\, A(x), \exists x\, A(x) \to A(a) \fCenter \Delta$
\end{prooftree}

By doing this, we obtain a proof of the desired sequent, together with many formulae of the form $A(a) \to \forall x\, A(x)$ or $\exists x\, A(x) \to A(a)$ on the left-hand side. However, we can eliminate each of them by adding an existential quantifier inference and cutting with formulae of the form
$$\vdash \exists y\, \big(A(y) \to \forall x\, A(x)\big)$$
or of the form
$$\vdash \exists y\, \big(\exists x\, A(x) \to A(y)\big),$$
both of which are easily derivable. Note that the existential quantifier inferences can be carried out in a way that is permissible by $\lk$ because the initial proof satisfied substitutability, weak regularity and the side-variable condition.
\endproof

\begin{cor}
If a sequent is derivable in $\calcb$, then it is already derivable in $\lk$.
\end{cor}

\section{Non-Elementary Speed-Up} \label{SectSpeedup}
Our strategy for proving Theorem \ref{main} is to show that $\calcb$ simulates a strong calculus that is already non-elementarily faster than $\lk$. It is related to \textit{quantifier shifts}, formulae such as:
$$(\forall x\, A \to B) \to \exists x\, (A \to B).$$

The main feature of the calculus is the addition of certain rules for quickly shifting quantifiers. 
Below, if $A$ is a subformula of $B$, we say that $B$ is a \textit{context} of $A$. We generally write $\kappa[A]$ to emphasize that $A$ occurs in a context $\kappa[\cdot]$; they should be thought as syntactical operators $A \mapsto \kappa[A]$.

\begin{defi}
The calculus $\lkqs$ is obtained by extending $\lk$ with the following rules:
\begin{align*}
\infer[]{\Gamma, \kappa[Q'x\, (A \lhd B)] \vdash \Delta}{\Gamma, \kappa[Qx\, A \lhd B] \vdash \Delta} &&
\infer[]{\Gamma, \kappa[Q'x\, (A \lhd B)] \vdash \Delta}{\Gamma, \kappa[A \lhd Qx\,B] \vdash \Delta}\\
\infer[]{\Gamma \vdash \Delta, \kappa[Q'x\, (A \lhd B)]}{\Gamma \vdash \Delta, \kappa[Qx\, A \lhd B]} &&
\infer[]{\Gamma \vdash \Delta, \kappa[Q'x\, (A \lhd B)]}{\Gamma \vdash \Delta, \kappa[A \lhd Qx\, B]}
\end{align*}
where $\kappa[\cdot]$ is a context, $\lhd \in \{\wedge, \vee, \to\}$ and $Q' = Q$, except if $\lhd$ is $\to$ and $Q$ is taken from the antecedent, in which case $Q'$ is opposite. We refer to these rules as \textit{deep quantifier shifts}. 
\end{defi}

Here, recall that, as per our conventions, the syntax of first-order logic includes separate symbols for free and bound variables.

\begin{proposition} \label{PropLKQS}
Cut-free $\calcb$ simulates cut-free $\lkqs$ double-exponentially, i.e., every $\lkqs$-provable sequent is $\calcb$-provable and there is a double exponential function that bounds the length of the least cut-free $\calcb$-proof of a $\calcb$-provable sequent in terms of its least cut-free $\lkqs$-proof.
\end{proposition}
\proof
Let $\pi$ be a cut-free $\lkqs$-proof of size $n$. We transform it into a cut-free $\calcb$-proof. Assume by induction that there is only one application of a quantifier shift and this is the last inference. We proceed by cases. For simplicity, we assume $\kappa$ does not change the polarity of the quantifiers. We also assume $\lhd$ is $\to$; the other connectives are treated similarly. \\

\noindent\textsc{Case I.} The last inference is: 
\begin{align*}
\infer[]{\Gamma, \kappa[\exists x\, (A(x) \to B)] \vdash \Delta}{\Gamma, \kappa[\forall x\, A(x) \to B] \vdash \Delta} 
\end{align*}
so that the proof has the following structure:
\begin{prooftree}
\def\fCenter{\mbox{\ $\vdash$\ }}
\AxiomC{$\vdots\sigma_2$}
\UnaryInf$\Gamma'' \fCenter \Delta'', A(a)$
\RightLabel{($*$)}
\UnaryInf$\Gamma'' \fCenter \Delta'', \forall x\, A(x)$
\UnaryInfC{$\vdots\sigma_1$}
\UnaryInf$\Gamma' \fCenter \Delta', \forall x\, A(x)$
\AxiomC{$\vdots\sigma_3$}
\UnaryInf$\Gamma', B\fCenter \Delta'$
\BinaryInf$\Gamma', \forall x\, A(x) \to B \fCenter \Delta'$
\UnaryInfC{$\vdots\sigma_0$}
\UnaryInf$\Gamma, \kappa[\forall x\, A(x) \to B] \fCenter \Delta$
\UnaryInf$\Gamma, \kappa[\exists x\, (A(x) \to B)] \fCenter \Delta$
\end{prooftree}
where the $\sigma_i$ denote subproofs. The subproof $\sigma_0$ might split into several branches not indicated in the diagram, each of which could potentially include its own copy of $\forall x\, A(x)$ being inferred and its own subproof $\sigma_3'$ with conclusion of the form $\Pi, B \vdash \Lambda$. Each branch can be dealt with the same way and so we will only focus on the indicated parts of $\pi$.

Similarly, the subproof $\sigma_1$ might split into several branches not indicated in the diagram, each of which could potentially include its own copy of $\forall x\, A(x)$ being inferred. This we need to bear in mind. Note that each of those copies of $\forall x\, A(x)$ is necessarily inferred from a different characteristic variable. In fact, the variable $a$ does not appear below $(*)$, by the regularity of $\lkqs$.

We modify the proof. Our approach is as follows: we would like to merge the subproofs $\sigma_1$ and $\sigma_2$ simply by postponing the inference $(*)$ as follows:
\begin{prooftree}
\def\fCenter{\mbox{\ $\vdash$\ }}
\AxiomC{$\vdots\sigma_2 + \sigma_1$}
\UnaryInf$\Gamma' \fCenter \Delta', A(a)$
\AxiomC{$\vdots\sigma_3$}
\UnaryInf$\Gamma', B\fCenter \Delta'$
\BinaryInf$\Gamma', A(a) \to B \fCenter \Delta'$
\RightLabel{($**$)}
\UnaryInf$\Gamma', \exists x\, ( A(x) \to B) \fCenter \Delta'$
\UnaryInfC{$\vdots\sigma_0$}
\UnaryInf$\Gamma, \kappa[\exists x\, (A(x) \to B)] \fCenter \Delta$
\end{prooftree}
The problem that might arise is that some occurrence of $\forall x\, A(x)$ that would be contracted in $\sigma_1$ with the indicated occurrence is unable to be contracted. We describe how to deal with each of them.

Notice that the problematic occurrence of $\forall x\, A(x)$ must originate from a strong-quantifier inference. Omit that inference and drag the unquantified formula $A(b)$ until after $(**)$, so that at that point we have a derivation of 
$$\Gamma'_0, A(b), \Gamma'_1, \exists x\, ( A(x) \to B), \vdash \Delta',$$
where $\Gamma' = \Gamma'_0, \Gamma'_1$. Add to that some exchanges to obtain:
$$\Gamma', \exists x\, ( A(x) \to B), A(b) \vdash \Delta',$$
and a subproof
\begin{align} \label{eqqsCaseI}
\infer[]{\Gamma', B \vdash \Delta'}{\vdots\sigma_3}
\end{align}
in order to infer 
$$\Gamma', \exists x\, ( A(x) \to B), A(b) \to B \vdash \Delta'.$$
with an application of $\exists L$ and a contraction, we are left again with 
$$\Gamma', \exists x\, ( A(x) \to B) \vdash \Delta',$$
to which we can apply $\sigma_0$ as desired. The proof grows by adding a copy of \eqref{eqqsCaseI}---of size at most $n$---for each problematic occurrence of $\forall x\, A(x)$ (of which there are at most $n$). Hence, it grows quadratically.\\

\noindent\textsc{Case II.} The last inference is:
\begin{align*}
\infer[]{\Gamma, \kappa[\forall x\, (A(x) \to B)] \vdash \Delta}{\Gamma, \kappa[\exists x\, A(x) \to B] \vdash \Delta} 
\end{align*}
The proof has the following structure:
\begin{prooftree}
\def\fCenter{\mbox{\ $\vdash$\ }}
\AxiomC{$\vdots\sigma_2$}
\UnaryInf$\Gamma'' \fCenter \Delta'', A(t)$
\UnaryInf$\Gamma'' \fCenter \Delta'', \exists x\, A(x)$
\UnaryInfC{$\vdots\sigma_1$}
\UnaryInf$\Gamma' \fCenter \Delta', \exists x\, A(x)$
\AxiomC{$\vdots\sigma_3$}
\UnaryInf$\Gamma', B\fCenter \Delta'$
\BinaryInf$\Gamma', \exists x\, A(x) \to B \fCenter \Delta'$
\UnaryInfC{$\vdots\sigma_0$}
\UnaryInf$\Gamma, \kappa[\exists x\, A(x) \to B] \fCenter \Delta$
\UnaryInf$\Gamma, \kappa[\forall x\, (A(x) \to B)] \fCenter \Delta$
\end{prooftree}
As in Case I, the subproof $\sigma_0$ could branch off into several subtrees each of which is taken care of in the same way, and so we only focus on the indicated parts of $\pi$.
We would like to merge the subproofs $\sigma_1$ and $\sigma_2$ as follows:
\begin{prooftree}
\def\fCenter{\mbox{\ $\vdash$\ }}
\AxiomC{$\vdots\sigma_2 + \sigma_1$}
\UnaryInf$\Gamma' \fCenter \Delta', A(t)$
\AxiomC{$\vdots\sigma_3$}
\UnaryInf$\Gamma', B\fCenter \Delta'$
\BinaryInf$\Gamma', A(t) \to B \fCenter \Delta'$
\UnaryInf$\Gamma', \forall x\, ( A(x) \to B) \fCenter \Delta'$
\UnaryInfC{$\vdots\sigma_0$}
\UnaryInf$\Gamma, \kappa[\forall x\, (A(x) \to B)] \fCenter \Delta$
\end{prooftree}
As before, we face the problem of circumventing a contraction of (possibly several occurrences of) $\exists x\, A(x)$ in $\sigma_1$ and solve it in the same way. The proof grows quadratically again. We note the following: it might happen that the term $t$ contains some free variable $a$. Hence, the proof of this case does not necessarily go through in $\lk$, as $a$ could be the characteristic variable of a strong-quantifier inference in $\sigma_1$. Note that the side-variable condition is satisfied. \\

\noindent\textsc{Case III.} The last inference is: 
\begin{align*}
\infer[]{\Gamma, \kappa\big[\forall x\, (A \to B(x))\big] \vdash \Delta}{\Gamma, \kappa\big[A \to \forall x\,B(x)\big] \vdash \Delta} 
\end{align*}
This is analogous to \textsc{Case II.}\\

\noindent\textsc{Case IV.} The last inference is: 
\begin{align*}
\infer[]{\Gamma, \kappa\big[\exists x\, (A \to B(x))\big] \vdash \Delta}{\Gamma, \kappa\big[A \to \exists x\,B(x)\big] \vdash \Delta} 
\end{align*}
This is analogous to \textsc{Case I.}\\

\noindent\textsc{Case V.} The last inference is: 
\begin{align*}
\infer[]{\Gamma \vdash \Delta, \kappa\big[\exists x\, (A(x) \to B)\big]}
{\Gamma \vdash \Delta, \kappa\big[\forall x\, A(x) \to B\big]} 
\end{align*}
The proof has the following form, modulo qualifications as in Cases I and II:
\begin{prooftree}
\def\fCenter{\mbox{\ $\vdash$\ }}
\AxiomC{$\vdots\sigma_2$}
\UnaryInf$\Gamma'', A(t) \fCenter \Delta''$
\RightLabel{($*$)}
\UnaryInf$\Gamma'', \forall x\, A(x) \fCenter \Delta''$
\UnaryInfC{$\vdots\sigma_1$}
\UnaryInf$\Gamma', \forall x\, A(x) \fCenter \Delta', B$
\UnaryInf$\Gamma' \fCenter \Delta', \forall x\, A(x) \to B$
\UnaryInfC{$\vdots\sigma_0$}
\UnaryInf$\Gamma \fCenter \Delta, \kappa\big[\forall x\, A(x) \to B\big]$
\UnaryInf$\Gamma \fCenter \Delta, \kappa\big[\exists x\, (A(x) \to B)\big]$
\end{prooftree}
We would like to postpone the inference $(*)$, thus merging $\sigma_1$ and $\sigma_2$:
\begin{prooftree}
\def\fCenter{\mbox{\ $\vdash$\ }}
\AxiomC{$\vdots\sigma_2 + \sigma_1$}
\UnaryInf$\Gamma', A(t) \fCenter \Delta', B$
\RightLabel{($**$)}
\UnaryInf$\Gamma' \fCenter \Delta', A(t) \to B$
\UnaryInf$\Gamma' \fCenter \Delta', \exists x\, ( A(x) \to B)$
\UnaryInfC{$\vdots\sigma_0$}
\UnaryInf$\Gamma \fCenter \Delta, \kappa[\exists x\, (A(x) \to B)]$
\end{prooftree}
However, we might have---similarly to the previous cases---an occurrence of $\forall x\, A(x)$ that would be contracted on the left-hand side as part of $\sigma_1$. For simplicity, assume there is only one. This occurrence is inferred from a formula $A(s)$. Postpone this occurrence until after $(**)$ so that at that point we have a derivation of 
$$\Gamma'_0, A(s), \Gamma'_1 \vdash \Delta', \exists x\, (A(x) \to B),$$
where $\Gamma' = \Gamma'_0, \Gamma'_1$. Add to that some exchanges to obtain:
$$\Gamma', A(s) \vdash \Delta', \exists x\, (A(x) \to B),$$
and a weakening:
\begin{align} \label{eqqsCaseV}
\infer[WR]{\Gamma', A(s) \vdash \Delta', \exists x\, (A(x) \to B), B}
{\Gamma', A(s) \vdash \Delta', \exists x\, (A(x) \to B)}
\end{align}
in order to infer 
$$\Gamma' \vdash \Delta', \exists x\, (A(x) \to B), A(s) \to B.$$
with an application of $\exists R$ and a contraction, we are left again with 
$$\Gamma' \vdash \Delta', \exists x\, (A(x) \to B),$$
to which we can apply the rest of $\sigma_0$ as desired. In this case the proof only grows linearly. As in \textsc{Case II}, the proof of this case would not necessarily work for $\lk$, but it is correct from the point of view of $\calcb$, as follows from the regularity of $\lkqs$ and the fact that the conclusion of the proof is $\pi$.\\

\noindent\textsc{Case VI.} The last inference is: 
\begin{align*}
\infer[]{\Gamma \vdash \Delta, \kappa\big[\forall x\, (A(x) \to B)\big]}
{\Gamma \vdash \Delta, \kappa\big[\exists x\, A(x) \to B\big]} 
\end{align*}
The proof has the following form, modulo qualifications as in Cases I and II:
\begin{prooftree}
\def\fCenter{\mbox{\ $\vdash$\ }}
\AxiomC{$\vdots\sigma_2$}
\UnaryInf$\Gamma'', A(a) \fCenter \Delta''$
\RightLabel{($*$)}
\UnaryInf$\Gamma'', \exists x\, A(x) \fCenter \Delta''$
\UnaryInfC{$\vdots\sigma_1$}
\UnaryInf$\Gamma', \exists x\, A(x) \fCenter \Delta', B$
\UnaryInf$\Gamma' \fCenter \Delta', \exists x\, A(x) \to B$
\UnaryInfC{$\vdots\sigma_0$}
\UnaryInf$\Gamma \fCenter \Delta, \kappa\big[\exists x\, A(x) \to B\big]$
\UnaryInf$\Gamma \fCenter \Delta, \kappa\big[\forall x\, (A(x) \to B)\big]$
\end{prooftree}
As before, we would like to postpone the inference $(*)$, thus merging $\sigma_1$ and $\sigma_2$:
\begin{prooftree}
\def\fCenter{\mbox{\ $\vdash$\ }}
\AxiomC{$\vdots\sigma_2 + \sigma_1$}
\UnaryInf$\Gamma', A(a) \fCenter \Delta', B$
\UnaryInf$\Gamma' \fCenter \Delta', A(a) \to B$
\UnaryInf$\Gamma' \fCenter \Delta', \forall x\, ( A(x) \to B)$
\UnaryInfC{$\vdots\sigma_0$}
\UnaryInf$\Gamma \fCenter \Delta, \kappa[\forall x\, (A(x) \to B)]$
\end{prooftree}
We deal with contractions of $\exists x\, A(x)$ on the left-hand side as in the previous case. The proof grows linearly.\\

\noindent\textsc{Case VII.} The last inference is: 
\begin{align*}
\infer[]{\Gamma \vdash \Delta, \kappa\big[\forall x\, (A \to B(x))\big]}
{\Gamma \vdash \Delta, \kappa\big[A \to \forall x\,B(x)\big]} 
\end{align*}
This is analogous to \textsc{Case VI.}\\

\noindent\textsc{Case VIII.} The last inference is: 
\begin{align*}
\infer[]{\Gamma \vdash \Delta, \kappa\big[\exists x\, (A \to B(x))\big]}{\Gamma \vdash \Delta, \kappa\big[A \to \exists x\,B(x)\big]} 
\end{align*}
This is analogous to \textsc{Case V.}\\

Hence, we have dealt with all the cases. Note that for each quantifier shifted, the proof grows at most quadratically. Since there are at most $n$ deep quantifier shifts, the size of the resulting proof is bounded by 
$$n^{2^n} = 2^{2^n\cdot \log(n)} \approx 2^{2^n}.$$
This finishes the proof.
\endproof

Theorem \ref{main} is a consequence of Proposition \ref{PropLKQS} and the following result:

\begin{theorem} \label{mainqs}
There is no elementary function bounding the length of the shortest cut-free $\lk$-proof of a formula in terms of its shortest cut-free $\lkqs$-proof.
\end{theorem}
\proof
We will make use of a very specific family of sequents $\{S_i\}_{i < \omega}$ described in \cite{baazleitsch94} and due to Statman \cite{statman79}, and specific $\lk$-proofs thereof. The sequents and the proofs themselves are not important for our proof. What is relevant is that they have the following properties:
\begin{enumerate}
\item the size of $S_i$ is polynomial in $i$;
\item there is no bound on the size of their smallest cut-free $\lk$-proofs that is elementary in $i$;
\item the size of these proofs (with cuts), however, is polynomially bounded in $i$;
\item all cut formulae are closed; we can also assume they are prenex by, e.g., \cite[Theorem 3.3]{baazleitsch94}.
\end{enumerate}
Let $\Gamma_i \vdash \Delta_i$ be one of the sequents. We modify the proof as follows: first, replace each cut
\begin{align*}
	\infer[Cut]{\Gamma \vdash \Delta}{\Gamma \vdash \Delta, A & A, \Gamma \vdash \Delta}
\end{align*}
with an application of $\to L$:
\begin{align*}
	\infer[\to L]{\Gamma, A\to A \vdash \Delta}{\Gamma \vdash \Delta, A & \Gamma, A \vdash \Delta}
\end{align*}
We are left with a cut free proof $\pi_0$ whose end sequent is of the form:
\begin{align} \label{eqSpeedup1}
A_0 \to A_0, \hdots, A_m \to A_m, \Gamma_i \vdash \Delta_i.
\end{align}

Choose any occurrence of a quantifier in $A_0$ that is not in the scope of another quantifier. Since it appears once in the antecedent of $A_0 \to A_0$ and once in the consequent, it appears once with each polarity. Because $A_0$ is assumed to be prenex, we can apply two deep quantifier shifts, so that we are left with a proof $\pi_0'$ of the sequent:
\begin{align} \label{eqSpeedup1.5}
\forall x^0_0\, \exists x^0_1\, \big(A_0' \to A_0'\big), \hdots, A_m \to A_m, \Gamma_i \vdash \Delta_i,
\end{align}
Continue choosing any occurrence of a quantifier in $A_0'$ that is not in the scope of another quantifier (within $A_0'$) and applying deep quantifier shifts in pairs until we obtain a proof $\pi_1$ of the sequent:
\begin{align*} 
\forall x^0_0\, \exists x^0_1\, \hdots \big(\hat A_0 \to \hat A_0\big), \hdots, A_m \to A_m, \Gamma_i \vdash \Delta_i,
\end{align*}
where $\hat A_0$ is quantifier-free and the prefix of $\hat A_0 \to \hat A_0$ consists of alternating quantifiers. Repeat this procedure for each of the $A_i$ to obtain a proof $\pi_m$ where each formula $A_i \to A_i$ is replaced by an expression of the form:
\begin{align*}
\forall x^i_0\, \exists x^i_1\, \hdots \big(\hat A_i \to \hat A_i\big),
\end{align*}
where $\hat A_i$ is quantifier free and the quantifier prefix is alternating.
Let $\hat\Gamma^i \vdash \hat\Delta^i$ be the sequent:
\begin{align} \label{eqSpeedup2}
\forall x^0_0\, \exists x^0_1\, \hdots \big(\hat A_0 \to \hat A_0\big), \hdots, \forall x^m_0\, \exists x^m_1\, \hdots \big(\hat A_m \to \hat A_m\big), \Gamma_i \vdash \Delta_i.
\end{align}
Note that, since the size of the initial proof was polynomial in $i$, there were polynomially many quantifiers in the proof; hence, we added only polynomially many deep quantifier shifts, and so the size of the resulting proof $\pi_m$ of \eqref{eqSpeedup2} is bounded polynomially in $i$. Moreover, it is cut-free. Consequently, it suffices to show:

\begin{claim} \label{ClaimSU}
There is no elementary function bounding the size of the smallest cut-free $\lk$-proofs of \eqref{eqSpeedup2}.
\end{claim} 
\proof
Let $\{\sigma_i\}_{i < \omega}$ be a sequence of such proofs. First, we transform it into a sequence of proofs of the Skolemizations of the sequents \eqref{eqSpeedup2}, so that for each sequent, each of the implications $\hat A_i \to \hat A_i$ remains only universally quantified. Then, by Herbrand's theorem, there are propositional proofs $\{\theta_i\}_{i < \omega}$ of the Herbrand sequents of \eqref{eqSpeedup2}, each of which is of the form:
\begin{align} \label{eqSpeedup3}
\bigwedge_{k}B^0_k, \hdots, \bigwedge_{k}B^m_k, \Gamma' \vdash \Delta'.
\end{align}
Moreover, the lengths of the Herbrand sequents are bounded exponentially by the lengths of $\{\theta_i\}_{i < \omega}$ (see \cite[Theorem 4.3]{baazleitsch94}). 
Each conjunct $B^p_k$ is a quantifier-free implication of the form:
\begin{align} \label{eqSpeedup4}
A(t_1, \hdots, t_l) \rightarrow A(s_1, \hdots, s_l),
\end{align}
for some terms $t_1, \hdots, t_l, s_1, \hdots, s_l$. The key point is that, in the process of Skolemizing the sequents \eqref{eqSpeedup2}, exactly one term in each pair $(t_j, s_j)$ was weakly quantified and the other (which was under the scope of the strong quantifier) was replaced by a Skolem function. Hence, either $t_j$ is of the form $f(s_j, \vec r)$ or $s_j$ is of the form $f(t_j, \vec r)$, for some tuple of terms $\vec r$.

Transform the sequent \eqref{eqSpeedup3} as follows: pick a pair $(t_j, s_j)$ of some disjunct $B^p_k$ such that the element thereof that was replaced by a Skolem term---say, $f(s_j, \vec r) = t_j$---does not have any other of $t_1, \hdots, t_l, s_1, \hdots, s_l$ as an argument. Such a term of course exists---it is the first term whose quantifier was  shifted in \eqref{eqSpeedup1.5}. Substitute $s_j$ for $t_j$ throughout the sequent. 

By repeating this process sufficiently-many times, each formula \eqref{eqSpeedup4} is transformed into a propositional tautology of the form $A(\vec t) \to A(\vec t)$. This means that the sequent
\begin{align*}
\vdash \bigwedge_{k}B^0_k, \hdots, \bigwedge_{k}B^m_k.
\end{align*}
is also a propositional tautology, whereby so too is
$\Gamma' \vdash \Delta'$,
from which follows that it is $\lk$-provable. By \cite[Theorem 2]{baazetal12}, there is an $\lk$-proof of the unskolemized sequent of length exponential in that of $\Gamma' \vdash \Delta'$. But this is impossible, since---by assumption---there exist no short proofs of the unskolemized sequent.
\endproof
\noindent This establishes Theorem \ref{mainqs} and thus Theorem \ref{main}.
\endproof

\section{Properties of the Calculi} \label{SectApplications}
We discuss some properties of the calculi augmented with unsound inferences. Some of the arguments are only sketched---we leave the details to the interested reader.

\subsection{$\calcint$ and $\ljcalc$}
Consider $\calcint$ and $\ljcalc$, the analogs of the extended calculi $\calcb$ and $\calc$ for intuitionistic logic. Specifically, these are the calculi obtained by restricting possible sequents to those with at most one formula on the right-hand side. The calculi are not sound for intuitionistic logic, as shown by any of the examples in the introduction. A consequence of this that they do not admit cut elimination. The argument makes use of an important consequence of cut elimination, \textit{subformula property}---every formula appearing in a cut-free proof of a sequent $S$ is a subformula of a formula in $S$.

\begin{proposition} \label{NoCutEliminationForLJ+}
$\calcint$ and $\ljcalc$ do not admit cut elimination.
\end{proposition}
\proof
Consider example \eqref{example2}. Since only one formula appears on the right-hand side of each sequent, this is an $\calcint$-proof. Suppose towards a contradiction that there were a cut-free proof of
\begin{equation} \label{eqcutelimination}
\vdash \exists x\, \big( A(x) \rightarrow A(f(x))\big).
\end{equation} 
By the subformula property, this derivation would consist entirely of subformulae of \eqref{eqcutelimination}. In particular, it would contain no strong quantifier inferences. Therefore, it would already be an $\lj$-proof. But this is impossible, as \eqref{eqcutelimination} is not intuitionistically valid.
\endproof

Proposition \ref{NoCutEliminationForLJ+} has the consequence that $\calcb$ and $\calc$ do not admit cut elimination by an algorithm that resembles Gentzen's. We give a rough definition of what we mean by this, but the reader may consult the appendix of \cite{baazhetzl11} for further details.

\begin{defi}
We say that a cut-eliminating procedure is \emph{Gentzen-style} if it is a transformation of proofs consisting of permutation of rules, substitution of free variables, reduction of cuts of a formula to cuts of its outermost subformulae, and absorption of axioms, i.e., elimination of cuts
\begin{align*}
\infer[Cut]{\Gamma \vdash \Delta, A}{\Gamma \vdash \Delta, A & A \vdash A}
\end{align*} 
by deleting the indicated occurrence of $A \vdash A$.
\end{defi}

The point in considering Gentzen-style cut elimination algorithms is that they transform intuitionistic proofs into cut-free intuitionistic proofs. Thus, we obtain:

\begin{cor}
$\calc$ and $\calcb$ do not admit any Gentzen-style cut elimination.
\end{cor}

As is well known, the following three formulae are the only quantifier shifts not derivable in $\lj$:
\begin{enumerate}
\item \label{QS1} $\forall x\, (A \vee B(x)) \vdash A \vee \forall x\, B(x)$;
\item \label{QS2}  $(\forall x\, A(x) \to B) \vdash \exists x\, (A(x) \to B)$;
\item \label{QS3} $(A \to \exists x\, B(x)) \vdash \exists x\, (A \to B(x))$.
\end{enumerate}

They give rise to a characterization of $\ljcalc$:
\begin{proposition}\label{ljchar}
A sequent is provable in $\ljcalc$ if and only if it is provable in  $\lj$ with all quantifier shifts added as axioms.
\end{proposition}
\proof
All three quantifier shifts \ref{QS1}--\ref{QS3} are provable in $\ljcalc$ (see, e.g., Example \ref{example25}). We show the converse.
Clearly, the sequents $\vdash \forall y\, A(y) \to \forall x\, A(x)$ and $\vdash \exists y\, A(y) \to \exists x\, A(x)$ are provable. Hence, using the quantifier shifts \ref{QS2} and \ref {QS3} one derives $\vdash \exists y\, (A(y) \to \forall x\, A(x))$ and $\vdash \exists y\, (\exists x\, A(x) \to A(y))$. 
Therefore, we can apply the proof of Theorem \ref{TheoremCorrectness} to obtain the desired result.
\endproof

\subsection{Skolemization} 
The following are applications of  (the proof of) Theorem \ref{main}:

\begin{proposition} \label{PropSKCuts}
If $\Gamma \vdash \Delta$ is $\calcb$-derivable (resp. $\calc$-derivable), then its Skolemization is $\calcb$-derivable (resp. $\calc$-derivable) in quadratically many steps using additional cuts.
\end{proposition}
\proof
Let $\Gamma \vdash \Delta$ be any $\calc$-derivable sequent. For each formula $A$ in $\Delta$, the sequent $A \vdash sk(A)$ is already $\lk$-derivable, it follows that, by using one cut per formula, we can derive $\Gamma \vdash sk(\Delta)$ in $\calc$. Similarly, $sk(A) \vdash A$ is already $\lk$-derivable. Hence, by also using one cut for each formula in $\Gamma$, we can derive $sk(\Gamma)\vdash sk(\Delta)$.
\endproof

However, as we see, the additional cuts are necessary:

\begin{proposition} \label{Propcor1}
Let $S$ be a sequent. Then, there is no elementary bound on the length of the smallest cut-free $\calc$-proof of the Skolemization of $S$ in terms of the smallest cut-free $\calcb$-proof of $S$.
\end{proposition}
\proof
We argue as in the proof of Claim \ref{ClaimSU}: if the Skolemization of $\Gamma \vdash \Delta$ were cut-free $\calc$-derivable in elementary many steps, the cut-free $\calc$-proof would already be an $\lk$-proof, whence by \cite[Theorem 2]{baazetal12}, there would be an $\lk$-proof of $\Gamma \vdash \Delta$ of elementary length with respect to the $\calcb$-proof.
\endproof
Note, in contrast, that a cut-free proof can be Skolemized in $\lk$ at no extra cost.

\begin{proposition}
There is no elementary bound on the length of the Herbrand sequent of a sequent in terms of its smallest cut-free $\calcb$-proof (resp. $\calc$-proof).
\end{proposition}
\proof
This follows from Proposition \ref{Propcor1}, as a cut-free $\calc$-proof of the Skolemization of a sequent contains no strong quantifier and is hence already an $\lk$-proof, whence we can apply Herbrand's theorem.
\endproof

Although deskolemization of cut-free proofs is exponential in $\lk$, it is linear in $\calc$:
\begin{proposition} \label{PropSKCuts2}
If the Skolemization of $\Gamma \vdash \Delta$ is cut-free $\calc$-derivable, then $\Gamma \vdash \Delta$ is cut-free $\calc$-derivable in linearly many steps.
\end{proposition}
\proof
Let $\pi$ be a cut-free $\calc$-proof of the Skolemization of $\Gamma \vdash \Delta$. By definition, the sets of Skolem functions assigned to each formula in $\Gamma \vdash \Delta$ are disjoint. Modify $\pi$ into a proof of $\Gamma \vdash \Delta$ by replacing each Skolem term $t$ by a new free variable $a_t$ and quantifying it as soon as it is introduced. This might not have been possible in $\lk$, but it is in $\calc$: substitutability clearly holds, as does weak regularity. The side-variable condition must also hold, since every Skolem function has only finite arity.
\endproof

\begin{remark}
The proof of Proposition \ref{PropSKCuts2} also yields the intuitionistic analog.
\end{remark} 

\section{Connections to the Epsilon Calculus} \label{SectEpsilon}
We start by recalling the sequent-calculus reformulation of Hilbert's \textit{$\varepsilon$-calculus} (see \cite{hilbertbernays}).
\begin{defi}
The quantifier-free language $\mathcal L_\varepsilon$ is obtained by removing symbols $\exists$ and $\forall$ from the language of first-order logic and adding symbols $\varepsilon$ and $\tau$. Formulae and terms are simultaneously defined by induction in a way that the following clauses are satisfied. We leave the precise recursion to the reader.
\begin{enumerate}
\item constants and free variables are terms;
\item if $t_1, \hdots, t_n$ are terms and $f$ is an $n$-ary function symbol, then $f(t_1, \hdots, t_n)$ is a term;
\item if $A(t)$ is a formula, $t$ is a term, and $x$ is a bound variable, then $\varepsilon_xA(x)$ and $\tau_xA(x)$ are terms.
\item if $t_0, \hdots, t_n$ are terms and $A$ is an $n$-ary predicate symbol, then $A(t_0, \hdots, t_n)$ is a formula;
\item if $A$ and $B$ are formulae, then $A\wedge B$, $A\vee B$, $\neg A$, and $A \to B$ are formulae.
\end{enumerate}
The \emph{$\varepsilon$-sequent calculus} $\lkep$ is obtained by adding to Propositional $\mathrm{LK}$ the rules:
\begin{align*}
	\infer[\tau]{\Gamma, A(\tau_x A(x)) \vdash \Delta}{\Gamma, A(t) \vdash \Delta} && 
	\infer[\varepsilon]{\Gamma \vdash \Delta, A(\varepsilon_x A(x))}{\Gamma \vdash \Delta, A(t)} 
\end{align*}
\end{defi}
The term $\varepsilon_x A$ is to be understood as a `generic' object satisfying property $A$. The dual term $\tau_x A$ is added for symmetry; it is in fact redundant---one can define $A(\tau_x A) = A(\varepsilon_x(\neg A))$. The \textit{standard translation} of first-order logic into epsilon calculus is defined recursively by mapping atomic formulae to themselves, respecting Booleans, and dealing with quantified formulae as follows: suppose $A'$ is the standard translation of $A$, then we map:
\begin{align*}
\exists x\, A(x) \mapsto A'(\varepsilon_xA'(x)), &&
\forall x\, A(x) \mapsto A'(\tau_xA'(x)).
\end{align*}

\begin{proposition} \label{PropEC}
A sequent is $\lk$-derivable if and only if the standard translation of its Skolemization is $\lkep$-derivable.
\end{proposition}

In the $\lkep$-proof of the standard translation of a sequent, $\varepsilon$- and $\tau$-rules play the role of weak quantifier inferences, while the substitution of variables for $\varepsilon$- and $\tau$-terms throughout the proof plays the role of strong quantifier inferences. The proof of Proposition \ref{PropEC} relies on the \textit{Extended First Epsilon Theorem} (see \cite[Theorem 16]{moserzach06}). It states the following: let
\begin{equation} \label{eqEFET}
A_0, \hdots, A_m \vdash A_{m+1}, \hdots, A_{m+n}
\end{equation}
be an $\lkep$-provable sequent all of whose occurrences of $\varepsilon$ and $\tau$ are \textit{weak} (this is defined analogously as for quantifiers in first-order logic) and let $\vec {t^i}$ be all the terms of the form $\tau_x A(x)$ or $\varepsilon_x A(x)$ in $A_i$. Then there exist finitely-many tuples of terms $\vec {t^i}_j$ without any occurrence of the symbols $\tau$ or $\varepsilon$ such that if $A^i_j$ is the result of substituting the terms $\vec {t^i}_j$ for $\vec {t^i}$ in $A_i$, then the Herbrand sequent
\begin{equation*}
\bigwedge_i A^i_0, \hdots, \bigwedge_i A^i_m \vdash \bigvee_i A^i_{m+1}, \hdots, \bigvee_i A^i_{m+n}
\end{equation*}
is a propositional tautology. 
If so, then the sequent resulting from `shifting' the disjunctions and conjunctions inside each formula to the subformula where the $\varepsilon$- and $\tau$-terms appeared originally is also a propositional tautology. 
Call this sequent $S$. Assuming \eqref{eqEFET} was the standard translation of the Skolemization of a sequent $\Gamma \vdash \Delta$, one can then use $S$ to obtain an $\lk$-proof of $\Gamma \vdash \Delta$. Let us illustrate this with an example:

\begin{exm}
Consider the following proof of a sequent whose only occurrence of $\tau$ is weak:
\begin{prooftree}
\def\fCenter{\mbox{\ $\vdash$\ }}
\Axiom$A(\tau_x(A(x)\to B)) \fCenter A(\tau_x(A(x)\to B))$
\UnaryInf$A(\tau_xA(x)) \fCenter A(\tau_x(A(x)\to B))$
\UnaryInf$A(\tau_xA(x)) \fCenter B,\, A(\tau_x(A(x)\to B))$
\UnaryInf$\fCenter A(\tau_xA(x)) \rightarrow B,\, A(\tau_x(A(x)\to B))$
\Axiom$B \fCenter B$
\BinaryInf$A(\tau_x(A(x)\to B)) \rightarrow B \fCenter A(\tau_xA(x)) \rightarrow B,\, B$
\UnaryInf$A(\tau_x(A(x)\to B)) \rightarrow B,\, A(\tau_xA(x)) \fCenter A(\tau_xA(x)) \rightarrow B,\, B$
\UnaryInf$A(\tau_x(A(x)\to B)) \rightarrow B \fCenter A(\tau_xA(x)) \rightarrow B,\, A(\tau_xA(x)) \rightarrow B$
\UnaryInf$A(\tau_x(A(x)\to B)) \rightarrow B \fCenter A(\tau_xA(x)) \rightarrow B$
\end{prooftree}

The corresponding Herbrand sequent is 
\begin{equation}\label{eqexampleepsilon}
A(a) \to B \vdash \big(A(a) \to B\big) \vee \big(A(b) \to B\big),
\end{equation}
which is a propositional tautology. Moreover, the result of shifting the disjunction to where the $\varepsilon$-term originally appeared, namely,
$$A(a) \to B \vdash A(a)\wedge A(b) \to B,$$
is also a propositional tautology. Here, the disjunction is replaced by a conjunction, as the $\varepsilon$-term appeared in the antecedent of an implication. The conclusion of the proof is the translation of the Skolemized sequent
$$\forall x\, \big(A(x) \to B\big) \vdash \forall x\, A(x) \to B,$$
and so this is $\lk$-provable by Proposition \ref{PropEC}. Compare this with
$$\forall x\,\big(A(x) \to B\big)\vdash \exists x\, \big(A(x) \to B\big),$$
which is the Skolemized sequent suggested by \eqref{eqexampleepsilon}.
\end{exm}

There is a very intimate connection between $\calc$ and the $\varepsilon$-calculus. In fact, the characteristic-variable conditions of $\calc$ are precisely chosen so as to ensure that we have the following result:

\begin{proposition} \label{propepsilon2}
If a sequent is has an $\calc$-proof of length $k$, then its standard translation has an $\lkep$-proof of length $\leq k$.
\end{proposition}

The proof proceeds by replacing weak-quantifier inferences  with $\varepsilon$- or $\tau$-rules and omitting strong quantifier inferences---instead, we substitute the appropriate $\varepsilon$- or $\tau$-term for the characteristic variable throughout the proof.

Proposition \ref{propepsilon2} readily yields an alternative proof of Theorem \ref{TheoremCorrectness}: let $\Gamma \vdash \Delta$ be an $\calc$-provable sequent. Then its Skolemization is $\calc$-provable by Proposition \ref{PropSKCuts}. But the Skolemization is common to $\lk$ and preserves validity, so by Propositions \ref{PropEC} and \ref{propepsilon2}, $\Gamma \vdash \Delta$ is $\lk$-provable. 

One can state the relationship between the $\varepsilon$-calculus and $\ljcalc$. 

\begin{proposition}
A sequent is $\ljcalc$-provable if and only if its standard translation is $\lj^{\varepsilon}$-provable.
\end{proposition}
\proof
Suppose $\Gamma \vdash D$ is $\ljcalc$-provable. By Proposition \ref{ljchar} it is then $\lj$-provable by adding quantifier shifts as axioms. As in Proposition \ref{propepsilon2}, the standard translation of $\Gamma \vdash D$ is $\lj^{\varepsilon}$-provable by adding the standard translations of quantifier shifts as axioms. However, the standard translations of quantifier shifts are already $\lj^{\varepsilon}$-provable.

Conversely, suppose the standard translation of $\Gamma \vdash D$ is $\lj^{\varepsilon}$-provable, say by a proof $\pi$. We want to directly translate $\pi$ into an $\ljcalc$-proof. The problem is that $\pi$ might have inferences arranged the following way:
\begin{prooftree}
\def\fCenter{\mbox{\ $\vdash$\ }}
\AxiomC{$\vdots$}
\UnaryInf$\Pi \fCenter B(s)$
\UnaryInf$\Pi \fCenter B(\varepsilon_x B(x))$
\UnaryInfC{$\vdots$}
\UnaryInf$\Pi' \fCenter A(\varepsilon_x B(x))$
\UnaryInf$\Pi' \fCenter A(\varepsilon_x A(x))$
\end{prooftree}

To rid ourselves of the problem, we will transform $\pi$ into a cut-free proof where no inferred $\varepsilon$- or $\tau$-term is modified by an $\varepsilon$- or $\tau$-inference. We do this by replacing occurrences of $A(\varepsilon_x A)$ by $A(f_A(x))$, for Skolem functions $f_A$. Suppose $A(t)$ is a cut formula in $\pi$ and $t$ is a term of the form $\varepsilon_x B$ that is inferred in $\pi$. For example:
\begin{prooftree}
\def\fCenter{\mbox{\ $\vdash$\ }}
\AxiomC{$\vdots$}
\UnaryInf$\Pi \fCenter B(s)$
\UnaryInf$\Pi \fCenter B(\varepsilon_x B(x))$
\UnaryInfC{$\vdots$}
\UnaryInf$\Pi' \fCenter A(\varepsilon_x B(x))$
\AxiomC{$\pi_0\ \vdots$}
\UnaryInf$A(\varepsilon_x B(x)), \Pi' \fCenter E'$
\BinaryInf$\Pi' \fCenter E'$
\UnaryInfC{$\vdots$}
\end{prooftree}
Modify $\pi_0$ into a proof $\pi_1$ as follows: the indicated occurrence of $\varepsilon_x B(x)$ must have been either introduced by weakening---in which case we modify the weakening so that $f_B(s)$ is introduced instead---or by an axiom:
\begin{prooftree}
\def\fCenter{\mbox{\ $\vdash$\ }}
\Axiom$C(\varepsilon_x B(x)) \fCenter C(\varepsilon_x B(x))$
\UnaryInfC{$\vdots$}
\UnaryInf$A(\varepsilon_x B(x)), \Pi' \fCenter E'$
\end{prooftree}
in which case we replace it by:
\begin{prooftree}
\def\fCenter{\mbox{\ $\vdash$\ }}
\Axiom$C(f_B(s)) \fCenter C(f_B(s))$
\Axiom$C(\varepsilon_x B(x)) \fCenter C(\varepsilon_x B(x))$
\BinaryInf$C(f_B(s)) \to C(\varepsilon_x B(x)), C(f_B(s)) \fCenter C(\varepsilon_x B(x))$
\UnaryInfC{$\vdots$}
\UnaryInf$A(f_B(s)), \Pi'' \fCenter E'$
\end{prooftree}
where $\Pi''$ has the form $C(f_B(s)) \to C(\varepsilon_x B(x)), \Pi'$. We modify $\pi$---say, in the latter case---as follows:
\begin{prooftree}
\def\fCenter{\mbox{\ $\vdash$\ }}
\AxiomC{$\vdots$}
\UnaryInf$\Pi \fCenter B(s)$
\AxiomC{$\vdots$}
\UnaryInf$B(f_B(s)) \fCenter B(f_B(s))$
\BinaryInf$B(s) \to B(f_B(s)), \Pi \fCenter B(f_B(s))$
\UnaryInfC{$\vdots$}
\UnaryInf$B(s) \to B(f_B(s)), \Pi' \fCenter A(f_B(s))$
\AxiomC{$\pi_1\ \vdots$}
\UnaryInf$A(f_B(s)), \Pi'' \fCenter E'$
\BinaryInf$C(f_B(s)) \to C(\varepsilon_x B(x)), B(s) \to B(f_B(s)), \Pi' \fCenter E'$
\UnaryInfC{$\vdots$}
\end{prooftree}

Repeat this for every weak $\varepsilon$- and $\tau$-term inferred in each cut formula. We then make a similar modification to the proof for each inferred $\varepsilon$- or $\tau$-term that is modified by an $\varepsilon$- or $\tau$-inference. For example,
\begin{prooftree}
\def\fCenter{\mbox{\ $\vdash$\ }}
\AxiomC{$\vdots$}
\UnaryInf$\Pi, \fCenter B(s)$
\UnaryInf$\Pi, \fCenter B(\varepsilon_x A(x))$
\UnaryInfC{$\vdots$}
\UnaryInf$\Pi', \fCenter A(\varepsilon_x B(x))$
\UnaryInf$\Pi', \fCenter A(\varepsilon_x A(x))$
\end{prooftree}
becomes
\begin{prooftree}
\def\fCenter{\mbox{\ $\vdash$\ }}
\AxiomC{$\vdots$}
\UnaryInf$\Pi, \fCenter B(s)$
\Axiom$B(f_B(s)) \fCenter B(f_B(s))$
\BinaryInf$B(s) \to B(f_B(s)), \Pi, \fCenter B(f_B(s))$
\UnaryInfC{$\vdots$}
\UnaryInf$B(s) \to B(f_B(s)), \Pi', \fCenter A(f_B(s))$
\UnaryInf$B(s) \to B(f_B(s)), \Pi', \fCenter A(\varepsilon_x A(x))$
\end{prooftree}

After this, replace all cuts with $\to L$ rules as in the proof of Theorem \ref{mainqs} and add a $\tau$-rule for each Skolem term and each $\varepsilon$- and $\tau$-term in each of the introduced implications. 

This yields a proof of a sequent $A_1', \hdots, A_n', \Gamma' \vdash D'$ such that:
\begin{enumerate}
\item $\Gamma'$ and $D'$ are respectively the standard translations of $\Gamma$ and $D$;
\item Each $A_i'$ is an implication of the form $A \to A$ with only $\tau$-terms, and so the standard translation of an implication of the form $A \to A$ preceded by a block of universal quantifiers (in particular, $A_1', \hdots, A_n', \Gamma' \vdash D'$ is the standard translation of a sequent $A_1, \hdots, A_n, \Gamma \vdash D$);
\item No inferred $\varepsilon$- or $\tau$-term is modified by an $\varepsilon$- or $\tau$-inference in the proof.
\end{enumerate} 
Using this, one can translate back the proof to an $\ljcalc$-proof of 
$$A_1, \hdots, A_n, \Gamma \vdash D$$ sequent by sequent. This is done by replacing $\varepsilon$- and $\tau$-rules with weak quantifier inferences and introducing strong quantifier inferences whenever a formula is the $\varepsilon$-translation of a strongly quantified formula. The point is that the $\varepsilon$- and $\tau$-terms that were not inferred in the proof can be replaced by free variables.
Finally, we obtain a proof of $\Gamma \to D$ by cutting the $A_i$.  We leave the details to the reader.
\endproof

\section{Concluding Remarks}\label{SectClosing}
Many related questions remain open. Perhaps one should more thoroughly investigate the the role that (individual) inferences (should) play in proofs, both from a philosophical and from a technical standpoint. We believe that it is not straightforward, as there are grounds for challenging the view of inferences simply as one-step subproofs of proofs, or---conversely---of proofs (only) as arbitrary concatenations of inferences.

Another problem is to explore the addition of different unsound rules to proof systems. Here, we gave two examples, but many more are possible. It is unknown exactly how faster the calculus $\calcb$ really is---the weakening of the usual characteristic-variable condition into substitutability, weak regularity, and the side-variable condition can perhaps be exploited further.

It is not clear whether $\calcb$ is any (or significantly) slower than $\calc$. It is also not clear in what precise relation they both stand to $\lkqs$ and $\lkep$. All these questions seem to promise fruitful lines of future research.\\

\subsection*{Acknowledgements} We would like to thank the anonymous referee for several comments that certainly improved the paper. This work was partially supported by FWF grants I-2671-N35 and W1255-N23.


\begin{thebibliography}{10}

\bibitem{baazhetzl11}
M. Baaz and S Hetzl,
\newblock On the Non-Confluence of Cut-Elimination,
\newblock {\em J. Symb. Logic}, 76 (2011), 313--340.

\bibitem{baazetal12}
M. Baaz, S Hetzl, and D. Weller,
\newblock On the Complexity of Proof Deskolemization,
\newblock {\em J. Symb. Logic}, 77 (2012), 669--686.

\bibitem{baazleitsch94}
M. Baaz and A. Leitsch,
\newblock On Skolemization and Proof Complexity,
\newblock {\em Fund. Inform.}, 20 (1994), 353--379.

\bibitem{baazleitsch00}
M. Baaz and A. Leitsch,
\newblock Cut Normal Forms and Proof Complexity,
\newblock {\em Ann. Pure Appl. Logic}, 97 (2000), 127--177.

\bibitem{hilbertbernays}
D. Hilbert and P. Bernays,
\newblock {\em Grundlagen der Mathematik II} (1970), 
\newblock Springer--Verlag,

\bibitem{moserzach06}
G. Moser and R. Zach,
\newblock The Epsilon Calculus and Herbrand Complexity,
\newblock {\em Stud. Logica}, 82 (2006), 133--155.

\bibitem{orevkov82}
V. P. Orevkov,
\newblock Lower Bounds for Increasing Complexity of Derivations after Cut Elimination {\em (in Russian)},
\newblock {\em J. Soviet Math.} (2006), 2337--2350.

\bibitem{statman79}
R. Statman,
\newblock Lower Bounds on Herbrand's Theorem,
\newblock {\em Proc. Amer. Math. Soc.} 75 (1979), 104--107.


\end{thebibliography}
\end{document}